\documentclass[10pt]{amsart}
\usepackage{a4}
\usepackage{amssymb}
\usepackage{amscd}
\newtheorem{thm}{Theorem}

\newtheorem{lem}[thm]{Lemma}
\newtheorem{prop}[thm]{Proposition}
\newtheorem{defn}[thm]{Definition}

\newtheorem{rem}[thm]{Remark}
\newtheorem*{tthm}{Theorem}

\numberwithin{thm}{section}
\numberwithin{equation}{section}

\newcommand{\norm}[1]{\left\Vert#1\right\Vert}
\newcommand{\abs}[1]{\left\vert#1\right\vert}

\newcommand{\n}{\mathbb{N}}

\newcommand{\Hi}{\mathcal{H}}

\begin{document}
\title{Weak type $(2,H)$ and weak cotype $(2,H)$ of operator spaces}
\author{Hun Hee Lee}
\address{Department of Mathematical Sciences, Seoul National University
          San56-1 Shinrim-dong Kwanak-gu Seoul 151-747, Korea}
\email{Lee.hunhee@gmail.com, bbking@amath.kaist.ac.kr}

\begin{abstract}
Recently an operator space version of type and cotype, namely type $(p,H)$ and cotype $(q,H)$ of operator spaces
for $1\leq p \leq 2\leq q \leq \infty$ and a subquadratic and homogeneous Hilbetian operator space $H$ were introduced and investigated by the author.
In this paper we define weak type $(2,H)$ (resp. weak cotype $(2,H)$) of operator spaces,
which lies strictly between type $(2,H)$ (resp. cotype $(2,H)$) and type $(p,H)$ for all $1\leq p <2$ (resp. cotype $(q,H)$ for all $2<q \leq \infty$).
This is an analogue of weak type 2 and weak cotype 2 in the Banach space case, so we develop analogous equivalent formulations.
We also consider weak-$H$ space, spaces with weak type $(2,H)$ and weak cotype $(2,H^*)$ simultaneously and establish
corresponding equivalent formulations.
\end{abstract}
\maketitle

\section{Introduction}
Hilbert spaces are usually regarded as best among Banach spaces in many aspects.
One of them is the characterization of Hilbert spaces by type and cotype notions.
Recall that a Banach space $X$ is called (gaussian) type $p$ ($1\leq p \leq 2$) if there is a constant $C>0$ such that
$$\pi_{q,2}(v^*) \leq C\cdot \ell^*(v)$$ for $\frac{1}{p} + \frac{1}{q} =1$ and any $v: X \rightarrow \ell^n_2$, $n\in \mathbb{N}$ 
and (gaussian) cotype $q$ if there is a constant $C'>0$ such that
$$\pi_{q,2}(u) \leq C'\cdot \ell(u)$$ for any $u: \ell^n_2 \rightarrow X$ and $n\in \mathbb{N}$.
Here, $\pi_{q,2}(\cdot)$ is the $(q,2)$-summing norm defined by
$$\pi_{q,2}(T:X\rightarrow Y) = \sup\Big\{ \frac{(\sum_k \norm{TSe_k}^q )^{\frac{1}{q}}}{\norm{S : \ell_2 \rightarrow X}} \Big\},$$ 
and $\ell^*(\cdot)$ is the trace dual of $\ell$-norm $\ell(\cdot)$ defined by
$$\ell(u) := \Big[\int_{\Omega} \norm{\sum^n_{k=1}g_k(\omega)ue_k}^2_X dP(\omega)\Big]^{\frac{1}{2}}$$
for a family of i.i.d. standard gaussian variables $\{g_k\}$ on a probability space $(\Omega, P)$.
It is well known that a Banach space is isomorphic to a Hilbert space if and only if it has type 2 and cotype 2 simultaneously.
Since type $p$ implies type $r$ for $1\leq r < p\leq 2$ and cotype $q$ implies cotype $s$ for $2\leq q < s \leq \infty$
we can say that type and cotype are measurements for being close to Hilbert spaces.

In the operator space category homogeneous Hilbertian operator spaces play a similar role of Hilbert spaces
in the Banach space category (See \cite{P1.7} and \cite{P3}). 
A Hilbertian operator space $H$ (i.e. $H$ is isometric to a Hilbert space) is called homogeneous if 
for every $u : H \rightarrow H$ we have $\norm{u}_{cb}  = \norm{u}$. Note that we usually assume that $H$ is infinite dimensional and separable. 
Thus, it is natural to expect type and cotype notions which could be appropriate measurements for being close to $H$.
In \cite{L} the author have considered type $(p,H)$ and cotype $(q,H)$ of operator spaces with additional assumption for $H$ to be subquadratic,
which are generalizations of OH-type 2 and OH-cotype 2 of Pisier in \cite{P1.2} (Precise definitions will be given in section \ref{sec-wtype2Hcotype2H},
and see \cite{GP, GP2, L} for other type and cotype notions of operator spaces).
We say that $H$ is subquadratic if for all orthogonal projections $\{P_i\}^n_{i=1}$ in $H$ with $I_H = P_1 + \cdots + P_n$ we have 
$$\norm{x}^2_{B(\ell_2)\otimes_{\min}H} \leq \sum^n_{i=1}\norm{I_{B(\ell_2)}\otimes P_i (x)}^2_{B(\ell_2)\otimes_{\min}H}$$
for any $x\in B(\ell_2)\otimes H$ (See p.82 of \cite{P2}). Note that $R_p = [R, C]_{\frac{1}{p}}$ ($1\leq p \leq \infty$)
(complex interpolation of the row and column Hilbert spaces) and $\min \ell_2$ (the minimal quantization (chapter 3 of \cite{P3}) of $\ell_2$)
are subquadratic, but $\max \ell_2$ (the maximal quantization of $\ell_2$) is not subquadratic (not even completely isomorphic to such a space).

As in the Banach space case we could say that type $(p,H)$ and cotype $(q,H)$ are measurements for being close to $H$ and $H^*$, respectively.
For example $H$ and $H^*$ have type $(2,H)$ and cotype $(2,H)$, respectively, and if we assume that $H^*$ is also subquadratic, then
an operator space has type $(2,H)$ and cotype $(2,H^*)$ if and only if it is completely isomorphic to $H$. See \cite{L} for the details.

In the Banach space case, there is an old question whether type $p$ for all $1\leq p <2$
(resp. cotype $q$ for all $2<q \leq \infty$) implies type 2 (resp. cotype 2). It is well known that this is not true.
We have an intermediate notion called weak type 2 (resp. weak cotype 2),
which lies strictly between type 2 (resp. cotype 2) and type $p$ for all $1\leq p <2$ (resp. cotype $q$ for all $2<q \leq \infty$).
We also have examples of Banach spaces with weak type 2 and weak cotype 2 simultaneously (i.e. weak Hilbert spaces)
which is not isomorphic to any Hilbert spaces.

In \cite{MP} and \cite{Pa} it is shown that weak type 2 is equivalent to the extendability of bounded linear maps
from a subspace of the original space into Hilbert spaces by applying a suitable projection with large rank
and weak cotype 2 is equivalent to the existence of Euclidean subspaces of proportional dimension in every finite dimensional subspace.
Precisely speaking, a Banach space $X$ is weak type 2 if and only if there are constants $0< \delta < 1$ and $C>0$ satisfying the following :
for any $S\subseteq E$ and any $u : S\rightarrow \ell^n_2$ there is a projection $P : \ell^n_2 \rightarrow \ell^n_2$ and a map
$v : X \rightarrow \ell^n_2$ such that $$\text{rk}P >\delta n, v|_S = Pu \,\, \text{and} \,\,\norm{v} \leq C\norm{u}$$
and weak cotype 2 if and only if there are constants $0 < \delta < 1$ and $C>0$ such that we can find $F_1 \subseteq F$ satisfying
$$d_{F_1} = d(F,\ell^{\text{dim}F}_2) \leq C \; \text{and} \; \text{dim}F_1 \geq \delta \text{dim}F$$
for any finite dimensional $F \subseteq X$, where $d(\cdot,\cdot)$ is the Banach-Mazur distance defined by
$$d(F,E) = \inf \{ \norm{u}\norm{u^{-1}} | u : E \rightarrow F,\; \text{isomorphism} \}.$$
Furthermore, there are many interesting equivalent formulations of weak Hilbert spaces (\cite{P0}).

Coming back to operator space setting it is also natural to consider operator space analogues of the above weak notions.
This is the theme of this paper, which is organized as follows.
In section \ref{sec-wtype2Hcotype2H} we will define weak type $(2,H)$ and weak cotype $(2,H)$.
From the definition weak type $(2,H)$ (resp. weak cotype $(2,H)$) clearly lies between cotype $(2,H)$
(resp. cotype $(2,H)$) and type $(p,H)$ for all $1 \leq p < 2$. (resp. cotype $(q,H)$ for all $2< q \leq \infty$).
In section \ref{sec-equi-cotype2} we prove that weak cotype $(2,H)$ is equivalent to the existence of
a subspace completely semi-isomorphic to $H$ of proportional dimension in every finite dimensional subspace.
We say that an operator space $E$ is completely semi-isomorphic to $H$ if
there is an completely bounded isomorphism $u : H(I) \rightarrow E$ with bounded inverse for some index set $I$.
Note that $H(I)$ is well-defined for any index set $I$ (\cite{P1.5}). We will simply write $H_n$ when $I = \{1, \cdots, n\}$.
More precisely, we will show the following theorem.
\begin{tthm}
$E$ has weak cotype $(2,H)$ if and only if there are constants $0 < \delta < 1$ and $C >0$ such that we can find $F_1 \subseteq F$ satisfying
$$d^{H^*}_{F_1,s} :=\inf\{\norm{w}_{cb}\norm{w^{-1}} | w : H^*_{\text{dim}F_1} \rightarrow F_1, \; \text{isomorphism}\} \leq C $$ and
$$\text{dim}F_1 \geq \delta \text{dim}F$$ for any finite dimensional subspace $F \subseteq E$.
\end{tthm}
In section \ref{sec-equi-type2} we consider the case of weak type $(2,H)$.
We prove that weak type $(2,H)$ is equivalent to the completely bounded (shortly, cb-) extendability of bounded linear maps
from a subspace of the second dual space into $H$ by applying a suitable projection with large rank.
More precisely, we will show the following theorem.
\begin{tthm}
$E$ has weak type $(2,H)$ if and only if there are constants $0< \delta < 1$ and $C>0$ satisfying
for any $S\subseteq E^{**}$ and any $u : S\rightarrow H_n$
there is a projection $P : H_n \rightarrow H_n$ of $\text{rank}>\delta n$ and $v : E^{**} \rightarrow H_n$ such that
$$v|_S = Pu\; \text{and}\; \norm{v}_{cb} \leq C\norm{u}.$$
\end{tthm}
In section \ref{sec-weak-H} we consider weak-$H$ spaces, spaces with weak type $(2,H)$ and weak cotype $(2,H^*)$ simultaneously
and establish 9 equivalent formulations including the following. Note that we need to assume that $H^*$ is also subquadratic.
\begin{tthm}
$E$ is a weak-$H$ space if and only if for any $0<\delta <1$ there is a constant $C >0$ with the following property :
for any finite dimensional $F\subseteq E$ we can find $F_1 \subseteq F$ and an onto projection $P :E \rightarrow F_1$ satisfying
$$d^H_{F_1,cb} := d_{cb}(F_1, H_{\text{dim}F_1}) \leq C, \; \text{dim}F_1 \geq \delta \text{dim}F\,\, \text{and} \,\, \norm{P}_{cb} \leq C.$$
\end{tthm}

Throughout this paper, we assume that the reader is familiar with basic concepts in Banach spaces (\cite{DJT, P1, TJ})
and operator spaces (\cite{ER,P3}).
In this paper $E$ and $H$ will be reserved for an operator space and
a fixed, infinite dimensional, separable, subquadratic and homogeneous Hilbertian operator space.
As usual, $B(E)$ and $CB(E)$ denote the set of all bounded linear maps and all cb-maps on an operator space $E$, respectively.

\section{Weak type $(2,H)$ and weak cotype $(2,H)$}\label{sec-wtype2Hcotype2H}

We start this section by reviewing the definition of type $(p,H)$ and cotype $(q,H)$.
We say that an operator space $E$ has type $(2,H)$ if there is a constant $C>0$ such that 
$$\pi^*_{2,H}(v^*) \leq C \cdot \ell^*(v)$$ for all $n\in \n$ and $v : E \rightarrow \ell^n_2$ and
cotype $(2,H)$ if there is a constant $C'>0$ such that $$\pi^*_{2,H}(u)\leq C' \cdot \ell(u)$$
for all $n\in \n$ and $v : \ell^n_2 \rightarrow E$, where $\pi^*_{2,H}(\cdot)$ is trace dual of $(2,H)$-summing norm $\pi_{2,H}(\cdot)$ defined by
$$\pi_{2,H}(T:E\rightarrow F) = \sup\Big\{ \frac{(\sum_k \norm{TSe_k}^2)^{\frac{1}{2}}}{\norm{S : H^* \rightarrow E}_{cb}} \Big\}.$$
Note that we have (\cite{L}) $$\pi^*_{2,H}(u) = \inf \norm{A}_{HS}\norm{B}_{cb},$$ where the infimum runs over all
possible factorization $u : \ell^n_2 \stackrel{A}{\longrightarrow} H^* \stackrel{B}{\longrightarrow} E$.

For $1\leq p <2$ we use approximation number approach to define type $(p,H)$ and cotype $(q,H)$.
We define the $k$-th cb-approximation number of $T : E \rightarrow F$ by $$a^o_k(T) := \inf \{ \norm{T-S}_{cb} : S \in CB(E,F), \text{rk}(S) < k\}.$$ 
(See \cite{O} for operator space versions of Gelfand and Kolmogorov numbers.)
We say that an operator space $E$ has type $(p,H)$ if there is a constant $C>0$ such that
$$\Big(\sum_k a^o_k(v)^q \Big)^{\frac{1}{q}} \leq C \cdot \ell^*(v)$$ for $\frac{1}{p} + \frac{1}{q} =1$ and $v : E\rightarrow H_n$, $n\in \n$, 
and cotype $(q,H)$ if there is a constant $C'>0$ such that
$$\Big(\sum_k a^o_k(u)^q \Big)^{\frac{1}{q}} \leq C' \cdot \ell(u)$$ for all $n\in \n$ and $u : H^*_n \rightarrow E$,
where $H_n$ is the $n$-dimensional subspace of $H$. We denote the infimums of such $C$ and $C'$ by $T_{p,H}(E)$ and $C_{q,H}(E)$, respectively.
The above definition is inspired by the classical equivalence between $(\sum_k a_k(u)^q)^{\frac{1}{q}}$ and $\pi_{q,2}(u)$ for any 
$u : \ell^n_2 \rightarrow X$, (Corollary 19.7 of \cite{TJ}) where $a_k(\cdot)$ is the $k$-th approximation number 
defined by $$a_k(u) = \inf \{ \norm{u-v} : v \in B(\ell_2,X), \text{rk}(v) < k\}.$$
Indeed, there is a constant $C>0$ such that 
$$\Big(\sum_k a_k(u)^q \Big)^{\frac{1}{q}} \leq \pi_{q,2}(u) \leq \frac{Cq}{q-2}\Big(\sum_k a_k(u)^q \Big)^{\frac{1}{q}}$$
for any $u : \ell^n_2 \rightarrow X$.

Now we define weak type $(2,H)$ and weak cotype $(2,H)$.

\begin{defn}\label{WeakOHtypecotype2-def}
\begin{itemize}
\item[(1)] $E$ is called \textbf{weak type $(2,H)$} if there is a constant
$C>0$ such that $$\sup_{k\geq 1} \sqrt{k}a^o_k(v) \leq C \cdot
\ell^*(v)$$ for all $n\in\mathbb{N}$ and $v : E \rightarrow H_n$.

\item[(2)] $E$ is called \textbf{weak cotype $(2,H)$} if there is a constant
$C>0$ such that $$\sup_{k\geq 1} \sqrt{k}a^o_k(u) \leq C' \cdot
\ell(u)$$ for all $n\in\mathbb{N}$ and $u : H^*_n \rightarrow E$.
\end{itemize}
The infimum of such $C$ and $C'$ will be denoted by $wT^H_2(E)$ and $wC^H_2(E)$, respectively.
\end{defn}

\begin{rem}{\rm
\begin{itemize}
\item[(1)]
If we replace all cb-norms in the definition of weak type $(2,H)$, weak cotype $(2,H)$
into the usual operator norms, then we recover the definition of weak type 2 and weak cotype 2 of Banach spaces,
which means that operator spaces with weak type $(2,H)$ and weak cotype $(2,H)$ also have weak type $2$ and weak cotype $2$ as Banach spaces,
respectively.
\item[(2)]
From the definition it is clear that weak type $(2,H)$ (resp. weak cotype $(2,H)$) lies between type $(2,H)$ (resp. cotype $(2,H)$)
and type $(p,H)$ for all $1 \leq p < 2$ (resp. cotype $(q,H)$ for all $2< q \leq \infty$).

Moreover, the above gap is strict. Let $X$ be the 2-convexified Tsirelson spaces by W. B. Johnson (and T. Figiel)
(See chapter 13 of \cite{P1} and \cite{CS}). It is well-known that $X$ is weak cotype 2 and type 2, but not cotype 2.
Now we consider $\min X$, and note that $$a^o_k(u) = a_k(u)\,\;\text{for any}\; u : H^*_n \rightarrow \min X.$$
Thus, $\min X$ is weak cotype $(2,H)$, but it is not cotype $(2,H)$ by Remark 3.5. of \cite{L}.
On the other hand, $\max X^* = (\min X)^*$ has weak type $(2,H)$ by Proposition \ref{duality} and the fact that $X$ is $K$-convex.
However, $\max X^*$ cannot have type $(2,H)$ since $X$ does not have cotype 2 as a Banach space.
See \cite{L2} for more examples by Tsirelson type constructions.
\end{itemize}
}
\end{rem}

We have the following duality between weak type $(2,H)$ and weak cotype $(2,H)$ as in the Banach space case.

\begin{prop}\label{duality}

\begin{itemize}
\item[(1)] $E$ has weak type $(2,H)$ if and only if $E^*$ is $K$-convex as a Banach space and weak cotype $(2,H)$.
Moreover, $$wC^H_2(E^*) \leq wT^H_2(E) \leq wC^H_2(E^*)K(E).$$

\item[(2)] $E^*$ has weak cotype $(2,H)$ if and only if $E$ is $K$-convex as a Banach space and weak cotype $(2,H)$.
\end{itemize}

\end{prop}
\begin{proof}
By Lemma 3.10 and succeeding remark of \cite{P1} we have $$\sup_{k\geq 1} \sqrt{k} a^o_k(v^*) \leq
wC^H_2(E^*)\ell(v^*) \leq wC^H_2(E^*)K(E)\ell^*(v)$$ for any $v : E \rightarrow H_n$, where $K(E)$
is the $K$-convexity constant of $E$ (see chapter 2 of \cite{P1}).
Thus, we have that $wT^H_2(E) \leq wC^H_2(E^*)K(E)$, and $wC^H_2(E^*) \leq wT^o_2(E)$ comes from the trivial fact 
$\ell^*(v) \leq \ell(v^*).$ This completes the proof of (1), and the proof for (2) is the same.
\end{proof}

\begin{rem}{\rm

When $H = OH$, the operator Hilbert space by Pisier, we have another version of type and cotype of operator spaces,
namely $S_p$-type and $S_q$-cotype defined as follows (\cite{L}). 
An operator space $E$ is said to have $S_p$-type ($1\leq p \leq 2$) if there is a constant $C >0$ 
such that $$\pi^o_{q,2}(v^*) \leq C \cdot \ell^*(v)$$ for $\frac{1}{p} + \frac{1}{q} =1$ and any $v: E \rightarrow OH_n$, $n\in \mathbb{N}$,
where $\pi^o_{q,2}(\cdot)$ is the completely $(q,2)$-summing norm defined by
$$\pi^o_{q,2}(T:E\rightarrow F) = \sup\Big\{ \frac{\norm{(TSe_{ij})}_{S_q(F)}}{\norm{S : S_2 \rightarrow E}_{cb}} \Big\}$$
and $S_q(F)$ is the vector valued Schatten class introduced in \cite{P2}.

$E$ is said to have $S_q$-cotype if there is a constant $C' >0$ such that
$$\pi^o_{q,2}(u) \leq C' \cdot \ell(u)$$ for every $n\in \mathbb{N}$ and $u: OH_n \rightarrow E$.
Note that $S_2$-type (resp. $S_2$-cotype) coincide with type $(2,OH)$ (resp. cotype $(2,OH)$)
and $S_p$-type (resp. $S_q$-cotype) implies type $(p, OH)$ (resp. cotype$(q,OH)$). 

In this particular case we have a stronger result, namely
weak type $(2,OH)$ and weak cotype $(2,OH)$) imply $S_p$-type for all $1\leq p <2$ and $S_q$-cotype for all $2< q \leq \infty$, respectively.

Indeed, it is enough to show that there is a constant $C_q > 0$ such that
$$\pi^o_{q,2}(u) \leq C_q \sup_{k\geq 1} \sqrt{k} a^o_k(u)$$ for all $n\in \mathbb{N}$ and $u : H^*_n \rightarrow E$. 
Now we assume that $\sup_{k\geq 1} \sqrt{k} a^o_k(u) =1 $ and choose $v_k : H^*_n \rightarrow E$
with $\text{rank} < 2^k$ and $\norm{u-v_k}_{cb} \leq a^o_{2^k}(u)$. Let $d_0 = v_o := 0$ and $d_k := v_k - v_{k-1}$,
so that $\text{rk}(d_k) < 2^k + 2^{k-1}$ and $$\norm{d_k}_{cb} \leq \norm{v_k -u}_{cb} + \norm{u - v_{k-1}}_{cb}
\leq 2^{-\frac{k}{2}} + 2^{-\frac{k-1}{2}}.$$ Note that $v_k = u$ if $2^k \geq n$ since $\text{rk}u \leq n$.
Thus, $u = \sum_{0 \leq k < 1+\log_2 n}d_k$ and by Lemma \ref{(q,2)-estimate}, which will be proved in the following, we have
\begin{align*}
\pi^o_{q,2}(u) & \leq \sum_{0 \leq k < 1+\log_2 n}\pi^o_{q,2}(d_k)
\leq \sum_{0 \leq k < 1+\log_2 n}\norm{d_k}_{cb}\text{rk} (d_k)^{\frac{1}{4} + \frac{1}{2q}}\\
& \leq 4\sum_{k}2^{\frac{1}{2}(\frac{1}{q} - \frac{1}{2})k} := C_q < \infty.
\end{align*}

}
\end{rem}

\begin{lem}\label{(q,2)-estimate}
For $2 \leq q \leq \infty$ and $n\in \mathbb{N}$ we have $$\pi^o_{q,2}(I_{OH_n}) \leq n^{\frac{1}{4} + \frac{1}{2q}}.$$
Furthermore, we have $$\pi^o_{q,2}(u) \leq \norm{u}_{cb}\text{rk} (u)^{\frac{1}{4} + \frac{1}{2q}}$$ for all $u : OH_n \rightarrow E.$
\end{lem}
\begin{proof} For the first statement we observe the following.
$$\pi^o_{q,2}(I_{OH_n}) = \norm{I_{2,q} \otimes I_{OH_n} : S_2 \otimes_{\min} OH_n \rightarrow S_q(OH_n)}_{cb},$$
where $I_{2,q} : S_2 \rightarrow S_q$ is the formal identity and $\otimes_{\min}$ is the minimal tensor product of operator spaces.
Here, we consider the natural operator space structure on Schatten classes from \cite{P2}.
Then, we have
\begin{equation}\label{(2,2)-estimate}
\pi^o_{2,2}(I_{OH_n}) = \pi_2(I_{\ell^n_2}) = \sqrt{n}
\end{equation}
and
\begin{equation}\label{(infinity,2)-estimate}
\begin{split}
& \pi^o_{\infty, 2}(I_{OH_n}) \leq n^{\frac{1}{4}}.
\end{split}
\end{equation}
Indeed, by Proposition 7.2 of \cite{P3} we have for any $\sum^n_{i = 1}x_i \otimes e_i \in S_2 \otimes OH_n$ that
\begin{align*}
\norm{\sum^n_{i = 1}x_i \otimes e_i}^2_{S_{\infty}\otimes_{\min}OH_n}
& = \norm{\sum^n_{i = 1}x_i \otimes \overline{x}_i}_{S_{\infty}\otimes_{\min}\overline{S_{\infty}}}
= \norm{\sum^n_{i = 1}x_i \otimes \overline{x}_i}_{B(\ell_2 \otimes_2 \overline{\ell}_2)}\\
& \leq \norm{\sum^n_{i = 1}x_i \otimes \overline{x}_i}_{S_2(\ell_2 \otimes_2 \overline{\ell}_2)}
= \norm{\sum^n_{i = 1}x_i \otimes \overline{x}_i}_{(S_2(\ell_2 \otimes_2 \overline{\ell}_2))^r},
\end{align*}
where $\overline{E}$ implies the conjugate operator space of $E$,
and $S_p(\Hi)$ and $\Hi^r$ imply the Schatten class and the row Hilbert space on a Hilbert space $\Hi$, respectively.
Since we have
\begin{align*}
(S_2(\ell_2 \otimes_2 \overline{\ell}_2))^r & = (S_2(\ell_2) \otimes_2 S_2(\overline{\ell}_2))^r = (S_2 \otimes_2 \overline{S_2})^r\\
& = (S_2)^r \otimes_{\min} (\overline{S_2})^r = (S_2)^r \otimes_{\min} \overline{(S_2)^r},
\end{align*}
we get
\begin{align*}
\norm{\sum^n_{i = 1}x_i \otimes e_i}_{S_{\infty}\otimes_{\min}OH_n} & \leq \norm{OH_n \rightarrow (S_2)^r ,\; e_i \mapsto x_i}_{cb}\\
& \leq \norm{OH_n \rightarrow S_2 , \; e_i \mapsto x_i}_{cb} \norm{J|_E}_{cb}\\
& = \norm{\sum^n_{i = 1}x_i \otimes e_i}_{S_2\otimes_{\min}OH_n}\cdot\norm{J|_E}_{cb},
\end{align*}
where $J : S_2 \rightarrow (S_2)^r$ is the formal identity and $E = \text{span}\{x_i\}\subseteq S_2$.
If we let $F =\text{ran}(J|_E) \subseteq (S_2)^r$, by homogeneity of $S_2$ and $(S_2)^r$ we have
$E \cong OH_n$ and $F \cong R_n$ completely isometrically.
Consequently, we have $$\norm{J|_E}_{cb} = \norm{J_n : OH_n \rightarrow R_n, \,\, e_i \mapsto e_i}_{cb} = n^{\frac{1}{4}},$$
which implies (\ref{(infinity,2)-estimate}).

By combining (\ref{(2,2)-estimate}) and (\ref{(infinity,2)-estimate}) by complex interpolation we get
$$\pi^o_{q,2}(I_{OH_n}) \leq n^{\frac{1}{4} + \frac{1}{2q}}.$$

Now we consider $u : OH_n \rightarrow E$ and let $P$ be the orthogonal projection on $OH_n$ onto $(\text{ker}u)^\perp$, then we have
\begin{align*}
\pi^o_{q,2}(u)
& \leq \sup\Big\{\frac{\norm{(uPSe_{ij})}_{S_q(E)}}{\norm{PS : S_2 \rightarrow (\text{ker}u)^\perp}_{cb}}\Big\}\\
& = \pi^o_{q,2}(u|_{(\text{ker}u)^\perp}I_{(\text{ker}u)^\perp}) \leq \norm{u}_{cb}\text{rk} (u)^{\frac{1}{4} + \frac{1}{2q}}
\end{align*}
by the first result.
\end{proof}

\section{Equivalent formulations of weak cotype $(2,H)$}\label{sec-equi-cotype2}

In this section we will establish equivalent formulations of weak cotype $(2,H)$.
We will see that weak cotype $(2,H)$ is equivalent to the existence of
a subspace completely semi-isomorphic to $H$ of proportional dimension in every finite dimensional subspace.
An operator space $E$ is said to be completely semi-isomorphic to an operator space $F$
if there is a completely bounded map $w : E \rightarrow F$ with bounded inverse and we denote 
$$d_s(E,F) = \inf\{\norm{w}_{cb}\norm{w^{-1}}\},$$ where the infimum runs over all such $w$
(See \cite{ORo} for the definition and other related topics). For $n\in \mathbb{N}$ and a $n$-dimensional operator space $E$
we consider $$d^{H^*}_{E,s} := d_s(H^*_n, E),$$ and we will see that $d^{H^*}_{E,s}$ is the right substitution for $d_E$ in this section
instead of $$d_{cb}(H^*_n, E) = \inf\{\norm{w}_{cb}\norm{w^{-1}}_{cb}| w : H^*_n \rightarrow E\}.$$

First, we start with an operator space version of Lemma 10.3 of \cite{P1} (or Lemma 2.2 of \cite{P0}).
The proof is almost the same with suitable modifications.
\begin{lem}\label{upgrading-estimate-lemma}
Let $0 < \alpha < 1$, $\beta > 0$ and $\lambda > 0$ be fixed constants and let $$\varphi(u) = \pi_{2,H}(u^*) \,\, \text{or} \,\,\ell(u).$$
Assume that for any $n\in \mathbb{N}$ and $u : H^*_n \rightarrow E$ we have
\begin{equation}\label{lem-assumption}
a^o_{[\alpha n]}(u) \leq \lambda [\alpha n]^{-\beta} \varphi(u)
\end{equation}
whenever $[\alpha n] \geq 1.$

Then there is a constant $C>0$ independent of the choice of $u$ and $\lambda$ such that we have
\begin{equation}\label{lem-conclusion}
\sup_{k\geq 1}k^{\beta}a^o_k(u) \leq \lambda C \cdot \varphi(u)
\end{equation}
for any $u : H^*_n \rightarrow E$.
\end{lem}
\begin{proof}
By homogeneity we may assume that $\lambda = 1$. Let $C_n$ be the best $C$ for which (\ref{lem-conclusion}) holds for all $u: H^*_n \rightarrow E$,
and we will show that $C_n$ is uniformly bounded for all $n\in \mathbb{N}$.

Now we fix $u : H^*_n \rightarrow E$ with $\varphi(u) \leq 1$, and let $N = [\frac{k}{\alpha}]$ for fixed $k\geq 1$.
Then there is a subspace $S \subseteq H^*_n$ such that
\begin{equation}\label{lem-pf1}
\norm{u|_S}_{cb} = a^o_{N+1}(u) \leq C_n(N+1)^{-\beta}
\end{equation}
and
\begin{equation}\label{lem-pf2}
a^o_{[\alpha N]}(u|_{S^\perp}) \leq [\alpha N]^{-\beta}.
\end{equation}

Indeed, when $N < n$, there is $v : H^*_n \rightarrow \text{ran}u \subseteq E$ with
$$\text{rk}v = N\; \text{and}\; \norm{u-v}_{cb} = a^o_{N+1}(u)$$ by a usual compactness argument.
Then $S = \text{ker}v$ is our desired subspace since $$a^o_{N+1}(u) \leq \norm{u - u|_{S^{\perp}}}_{cb} = \norm{u|_S}_{cb}
= \norm{(u-v)|_S}_{cb} \leq \norm{u-v}_{cb} = a^o_{N+1}(u)$$ and by (\ref{lem-assumption}) we have
$a^o_{[\alpha N]}(u|_{S^\perp}) \leq [\alpha N]^{-\beta}.$

When $N \geq n$, $S = \{0\}$ is our desired subspace since $\norm{u|_S}_{cb} = a^o_{N+1}(u) = 0$ and by (\ref{lem-assumption}) we have
$$a^o_{[\alpha N]}(u|_{S^\perp}) = a^o_{[\alpha N]}(u) = a^o_{[\alpha N]}(uPi) \leq a^o_{[\alpha N]}(uP)
\leq [\alpha N]^{-\beta}\varphi(uP) \leq [\alpha N]^{-\beta},$$
where $i : H^*_n \hookrightarrow H^*_N$ and $P$ is the orthogonal projection from $H^*_N$ onto $H^*_n$.

Thus, by combining (\ref{lem-pf1}) and  (\ref{lem-pf2}) we have
\begin{align*}\label{pf-estimate}
a^o_{[\alpha N]}(u) & \leq \inf\{ \norm{u-v}_{cb} : \text{rk}(v) < [\alpha N] \}\\ 
& \leq \inf\{ \norm{(u-v)|_S}_{cb} + \norm{(u-v)|_{S^\perp}}_{cb} : \text{rk}(v) < [\alpha N] \}\\
& \leq \norm{u|_S}_{cb} + \inf\{ \norm{u|_{S^\perp}- v|_{S^\perp}}_{cb} : \text{rk}(v) < [\alpha N], v|_S \equiv 0 \}\\
& \leq C_n(N+1)^{-\beta} + [\alpha N]^{-\beta}
\end{align*}
whenever $[\alpha N] \geq 1$, which is guaranteed if $k \geq 2.$

Since $k \geq [\alpha N]$, we have $$a^o_k(u) \leq a^o_{[\alpha N]}(u),$$ and $N+1 \geq \frac{k}{\alpha}$ implies 
$[\alpha N] \geq k - \alpha - 1$, so that $$k^{\beta}a^o_k(u) \leq C_n\alpha^{\beta} + (k(k - \alpha - 1)^{-1})^{\beta} 
\leq C_n\alpha^{\beta} + 3^{\beta}$$ for $k \geq 2$.

Since the case $k = 1$ is trivial by taking $n_0 \in \mathbb{N}$
such that $[\alpha n_0] = 1$ in (\ref{lem-assumption}), we get
$\sup_{k\geq 1} k^{\beta}a^o(u) \leq C_n\alpha^{\beta} + 3^{\beta}$ and consequently
$$C_n \leq C_n\alpha^{\beta} + 3^{\beta},$$ which lead us to the desired conclusion.

\end{proof}

The following is a partial relationship between $\pi^*_{2,H}$-norm and $\ell$-norm.

\begin{lem}\label{com-2-sum-ell}
Let $n\in \n$ and $u : H^*_n \rightarrow E$. Then we have $$\pi^*_{2,H}(u) \leq d^{H^*}_{E,s} \ell(u).$$
\end{lem}
\begin{proof}
Let $\text{dim}E = m \in \n$ and $T : E \rightarrow H^*_m$ be an arbitrary isomorphism.
Then we have
\begin{align*}
\pi^*_{2,H}(u) &\leq \norm{T^{-1}}_{cb}\pi^*_{2,H}(Tu) = \norm{T^{-1}}_{cb}\ell(Tu)\\
& \leq \norm{T^{-1}}_{cb}\norm{T} \ell(u).
\end{align*}
Since we take $T$ arbitrarily we get the desired result.
\end{proof}

Now we are ready to prove our main result in this section.
We mainly follow the proof of Theorem 10.2 of \cite{P1} with a slight modification.

\begin{thm}\label{locweakOHcotype2} The followings are equivalent.

\begin{itemize}
\item[(1)] $E$ has weak cotype $(2,H)$.

\item[(2)] There are constants $0 < \delta < 1$ and $C >0$ such that we can find $F_1 \subseteq F$ satisfying
$$d^{H^*}_{F_1,s} \leq C \;\text{and}\; \text{dim}F_1 \geq \delta \text{dim}F$$ for any finite dimensional subspace $F \subseteq E$.

\item[(3)] For all $0 < \delta < 1$, there is a constant $C> 0$ satisfying the same statement as in (2).
\end{itemize}
\end{thm}
\begin{proof}

(3) $\Rightarrow$ (2) : Trivial.

\vspace{0.3cm}

(2) $\Rightarrow$ (1) : Let $u : H^*_n \rightarrow E$ and $\alpha = 1 -\frac{\delta}{2}$. Now we claim that
\begin{equation}\label{claim1}
a^o_{[\alpha n]}(u) \leq C'\ell(u) [\alpha n]^{-\frac{1}{2}}
\end{equation}
for some $C'>0$ whenever $[\alpha n] \geq 1$. Note that we may assume that $u$ is invertible by a typical perturbation argument.
Let $F= \text{ran}u$. Then there is $F_1 \subseteq F$ with $$\text{dim}F_1 \geq \delta n\; \text{and} \; d^{H^*}_{F_1,s} \leq C.$$
If we set $S_1 = u^{-1}(F_1)$, then $u|_{S_1} : S_1 \rightarrow F_1$ satisfies the following by Lemma \ref{com-2-sum-ell}.
$$\Big(\sum_{k \geq 1}a^o_k(u|_{S_1})^2 \Big)^{\frac{1}{2}} \leq C\cdot \ell(u|_{S_1}) \leq C \cdot \ell(u),$$ which implies
$$a^o_k(u|_{S_1}) \leq C \cdot \ell(u)k^{-\frac{1}{2}}$$ for all $k\geq 1$.
If we let $N = n -\text{dim}S_1$, then the above is equivalent to $$a^o_{N+k}(u) \leq C \cdot \ell(u)k^{-\frac{1}{2}}$$ for all $k\geq 1$.
Choosing $k = [\frac{\delta n}{2}]$ we find $N + k \leq [\alpha n]$. Thus for $\frac{\delta n}{2} > 1$ we have
$$a^o_{[\alpha n]}(u) \leq C \cdot \ell(u) \Big[\frac{\delta n}{2}\Big]^{-\frac{1}{2}}
\leq C \cdot \ell(u) \Big(\frac{\delta n}{2} - 1 \Big)^{-\frac{1}{2}}.$$
Since $n \geq \frac{1}{\alpha} = \frac{2}{2-\delta}$, we have
$$\frac{\delta^2}{2-\delta} \leq \frac{\delta(1-\alpha)n}{4},$$ which
implies $$\frac{\delta n}{2} - 1 \geq \frac{\delta^3}{2-\delta}\alpha n \geq \frac{\delta^3}{2-\delta}[\alpha n].$$

When $\frac{\delta n}{2} \leq 1$, we have $n \leq \frac{2}{\delta}$, so that
$$a^o_{[\alpha n]}(u) \leq \norm{u}_{cb} \leq d_{cb}(H^*_n, \max{\ell_2^n})\norm{u} \leq n\ell(u) \leq \frac{2}{\delta}\ell(u).$$
Combining these we get (\ref{claim1}) for $C' = \max\{C \sqrt{\frac{2-\delta}{\delta^3}}, \frac{2}{\delta}\}$.

Now by (\ref{claim1}) and Lemma \ref{upgrading-estimate-lemma} we get $$\sup_{k\geq 1}\sqrt{k}a^o_k(u) \leq C'\ell(u)$$
for any $n \in \mathbb{N}$ and $u : H^*_n \rightarrow E$, which implies $E$ has weak cotype $(2,H)$.

\vspace{0.3cm}
(1) $\Rightarrow$ (3) :

Consider $F \subseteq E$ with $\text{dim}F = n$. Then by Theorem 3.11 of \cite{P1} there is an isomorphism $u : H_n^* \rightarrow F$
with $$\ell(u) = \sqrt{n} \,\,\text{and}\,\, \ell((u^{-1})^*)\leq K(1+\log d_F)\sqrt{n}$$ for some $K>0$.

Fix $1 \leq k \leq n$. Then by Theorem 5.8 of \cite{P1} there is $F_1 \subseteq F$ with $\text{codim}F_1 = k-1$ such that
$$\norm{u^{-1}|_{F_1}} \leq K(1+\log d_F) \sqrt{\frac{n}{k}}.$$
Consider $S_1 = u^{-1}(F_1)$ and $u|_{S_1} : S_1 \rightarrow F$. Then there is $S_2 \subseteq S_1$ with $$\text{dim}S_1/S_2 = k-1$$ such that
\begin{align*}
\norm{u|_{S_2}}_{cb} & \leq 2a^o_k(u|_{S_1}) \leq 2wC^H_2(E)\ell(u|_{S_1})k^{-\frac{1}{2}}\\
& \leq 2wC^H_2(E)\ell(u)k^{-\frac{1}{2}} = 2wC^H_2(E)\sqrt{\frac{n}{k}}.
\end{align*}
Now we set $F_2 = u(S_2)\subseteq F_1$. Then we have
$$\text{dim}F_2 = \text{dim}F_1 - \text{dim}(\text{ker}u|_{S_2}) = n-k+1-k+1 = n-2k+2.$$ Thus, we have
\begin{align*}
d^{H^*}_{F_2,s} & \leq \norm{u|_{S_2}}_{cb}\norm{u^{-1}|_{F_2}}\\
& \leq 2wC^H_2(E)\sqrt{\frac{n}{k}}K(1+\log d_F) \sqrt{\frac{n}{k}}\\
& \leq 2KwC^H_2(E)\Big(\frac{n}{k}\Big)(1+\log d^{H^*}_{F,s}).
\end{align*}
By replacing $2k-2$ into $k$ and observing $\frac{n}{k} \leq 2\cdot \frac{n}{2k-2}$ we get $\widetilde{F} \subseteq F$ with
$$\text{dim}\widetilde{F} = n-k$$ such that $$d^{H^*}_{\widetilde{F},s} \leq K'wC^H_2(E)\Big(\frac{n}{k}\Big)(1+\log d^{H^*}_{F,s})$$
for each $k =1,2,\cdots, n$ and some $K'>0$.

Now we define $$f(\delta) := \inf\{ d^{H^*}_{F_1,s} : F_1 \subseteq F, \text{dim}F_1 \geq \delta n \}$$ for $0 < \delta <1$.
Note that it is enough to show that $f(\delta) \leq B(\delta)$ for some constant $B(\delta) > 0$ depending only on $\delta.$

Consider any $F_1 \subseteq F$ with $\text{dim}F_1 = m \geq \delta n$. If $\delta m \leq m -1$, then we can choose
$1 \leq k_0 \leq m$ such that $$m - k_0 - 1 < \delta m \leq m - k_0$$ or equivalently
$$\Big( \frac{m}{2k_0} \leq \Big) \frac{m}{k_0 + 1} < (1-\delta)^{-1} \leq \frac{m}{k}.$$
Then by the preceding result we have $\widetilde{F} \subseteq F_1$ with $\text{dim}\widetilde{F} = m - k_0$ such that
\begin{align*}
d^{H^*}_{\widetilde{F},s} & \leq K'wC^H_2(E)\Big(\frac{m}{k_0}\Big)(1+\log d^{H^*}_{F_1,s})\\
& \leq 2K'wC^H_2(E)(1 - \delta)^{-1}(1+\log d^{H^*}_{F_1,s})
\end{align*}
and $\text{dim}\widetilde{F} = m - k_0 \geq \delta m \geq \delta^2 n$.

When $\delta m > m -1$ (equivalently $(1-\delta)^{-1} > m$), we set $\widetilde{F} = F_1$. Then we have
\begin{align*}
d^{H^*}_{\widetilde{F},s} & = d^{H^*}_{F_1,s} \leq m \leq m(1+\log d^{H^*}_{F_1,s})\\
& \leq (1 - \delta)^{-1}(1+\log d^{H^*}_{F_1,s}),
\end{align*}
and $\text{dim}\widetilde{F} = m > \delta m \geq \delta^2 n$.

By taking infimum over such $F_1$ we get $$f(\delta^2) \leq A(1 - \delta)^{-1}(1+\log f(\delta)),$$ where $A = \max\{1, 2K'wC^H_2(E)\}$.
Since we have $1\leq f(\delta) \leq \sqrt{n}$ trivially, we can apply Lemma 8.6 of \cite{P1}. Thus, we get
$$f(\delta) \leq K''A(1 - \delta)^{-1}(1+\log (K''A(1-\delta)^{-1})) := B(\delta)$$ for some $K''>0$.

\end{proof}

\section{Equivalent formulations of weak type $(2,H)$}\label{sec-equi-type2}

In this section we provide equivalent formulations of weak cotype $(2,H)$. We will see that weak type $(2,H)$ is equivalent to
the completely bounded extendability of bounded operators from a subspace of the second dual space into $H$
by applying a suitable projection with large rank. First, we establish some implications of weak type $(2,H)$.
In proving (1) $\Rightarrow$ (2) $\Rightarrow$ (3) of Theorem \ref{locweakOHtype2} we follow the proof of Theorem 11.6 in \cite{P1}.
However, we can not imitate the proof of (iv) $\Rightarrow$ (i) directly, in order to establish (3) $\Rightarrow$ (1),
due to the absence of appropriate operator space version of quotient of subspace (shortly QS) Theorem (See Theorem 8.4 of \cite{P1}).
Thus, we take another route using the formulation in the style of Theorem 3.1 in \cite{P0}.
Before that we prepare lemmas similar to Lemma \ref{upgrading-estimate-lemma},
and we omit the proofs since they are almost the same as before.

\begin{lem}\label{upgrading-estimate-lemma2}
Let $0 < \alpha < 1$, $\beta > 0$ and $\lambda > 0$ be fixed constants. Assume that for any $n\in \mathbb{N}$,
$u : \ell^n_2 \rightarrow E$ and $v : E \rightarrow \ell^n_2$ we have
$$a_{[\alpha n]}(vu) \leq \lambda [\alpha n]^{-\beta}\ell^*(v)\pi_{2,H}(u^*)$$ whenever $[\alpha n] \geq 1.$

Then there is a constant $C>0$ independent of the choice of $u$, $v$ and $\lambda$ such that we have
$$\sup_{k\geq 1}k^{\beta}a_k(vu) \leq \lambda C \ell^*(v)\pi_{2,H}(u^*)$$ for any $u$ and $v$.
\end{lem}

\begin{lem}\label{upgrading-estimate-lemma3}
Let $0 < \alpha < 1$, $\beta > 0$ and $\lambda > 0$ be fixed constants. Assume that for any $n\in \mathbb{N}$
and $v : E \rightarrow H_n$ we have
$$a^o_{[\alpha n]}(v) \leq \lambda [\alpha n]^{-\beta}\ell^*(v)$$ whenever $[\alpha n] \geq 1.$

Then there is a constant $C>0$ independent of the choice of $u$ and $\lambda$ such that we have
$$\sup_{k\geq 1}k^{\beta}a^o_k(v) \leq \lambda C \ell^*(v)$$ for any $v$.
\end{lem}

We also need the following operator space version of Weyl numbers and their basic properties.

\begin{defn}\label{CB-Weyl}
Let $T : E \rightarrow F$ be a cb-map between operator spaces.
Then we define the k-th \textbf{H-Weyl number} of $T$ by
$$x^H_k(T) := \sup \{ a^o_k(TA) : A \in CB(H^*, E), \norm{A}_{cb} \leq 1\}$$ for $k = 1, 2, \cdots.$
\end{defn}

\begin{prop}\label{Prop-Weyl-multiplicative}
Let $T : E \rightarrow F$ and $S : F \rightarrow G$ be cb-maps between operator spaces. Then we have
$$x^H_{m+n-1}(ST) \leq x^H_m(T)x^H_n(S)$$ for $n,m = 1, 2, \cdots.$
\end{prop}
\begin{proof}
Consider $A \in CB(H^* ,E)$ with $\norm{A}_{cb} \leq 1$. For any $\epsilon > 0$ we have $A_1 \in CB(H^* ,E)$ and $A_2 \in CB(H^* ,F)$
with $\text{rk}A_1 < m$ and $\text{rk}A_2 < n$ such that $$\norm{TA - A_1}_{cb} \leq (1+\epsilon)a^o_m(TA)$$
and $$\norm{S(TA - A_1)-A_2}_{cb} \leq (1+\epsilon)a^o_n(S(TA - A_1)).$$ Since $\text{rk}(SA_1 + A_2) < m+n-1$, we have
\begin{align}\label{cb-approx-multiplicative}
a^o_{m+n-1}(STA) & \leq \norm{STA - SA_1 - A_2}_{cb} \leq (1+\epsilon)a^o_n(S(TA - A_1))\nonumber \\
& \leq (1+\epsilon)\norm{TA - A_1}_{cb}x^o_n(S) \leq (1+\epsilon)^2 a^o_m(TA)x^o_n(S)\\
& \leq (1+\epsilon)^2 x^H_m(T)x^H_n(S).\nonumber
\end{align}
Since we took $\epsilon > 0$ arbitrarily, we get our desired result.
\end{proof}

\begin{prop}\label{Prop-Weyl}
Let $T : E \rightarrow H^*$ be a cb-map. Then we have $$\sqrt{k}x^H_k(T) \leq \pi_{2,H}(T)$$ for $k = 1, 2, \cdots.$
\end{prop}
\begin{proof}
Consider $A : H^* \rightarrow E$ with $\norm{A}_{cb} = \norm{\sum_i Ae_i \otimes e_i}_{E\otimes_{\min} H} \leq 1$.
Then for any $k = 1, 2, \cdots$ we have
\begin{align*}
\pi_{2,H}(T)& \geq \frac{(\sum_i \norm{TAe_i}^2)^{\frac{1}{2}}} {\norm{\sum_i Ae_i \otimes e_i}_{E\otimes_{\min} H}}\\
& \geq \Big(\sum_i \norm{TAe_i}^2\Big)^{\frac{1}{2}} = \norm{TA}_{HS} = \Big(\sum_i a^o_i(TA)^2 \Big)^{\frac{1}{2}}\\
& \geq \sqrt{k}a^o_k(TA),
\end{align*}
where $\norm{\cdot}_{HS}$ refers to the Hilbert-Schmidt norm.
\end{proof}

\begin{thm}\label{locweakOHtype2}
We have (1) $\Rightarrow$ (2) $\Rightarrow$ (3) $\Rightarrow$ (4).
\begin{itemize}
\item[(1)] $E$ has weak type $(2,H)$.

\item[(2)] There is a constant $C>0$ satisfying for any subspace $S\subseteq E$ and any $u : S\rightarrow H_n$
there is an extension $\widetilde{u} : E \rightarrow H_n$ such that for any $k=1,2,\cdots,n$ there is a projection
$P_k : H_n \rightarrow H_n$ with $\text{rank} > n - k$ and $$\norm{P_k\widetilde{u}}_{cb} \leq C\sqrt{\frac{n}{k}}\norm{u}.$$

\item[(3)] There are constants $0< \delta < 1$ and $C_{\delta}>0$ satisfying for any $S\subseteq E$ and any $u : S\rightarrow H_n$
there is a projection $P : H_n \rightarrow H_n$ of $\text{rank}>\delta n$ and $v : E \rightarrow H_n$ such that
$$v|_S = Pu \,\, \text{and} \,\, \norm{v}_{cb} \leq C_{\delta}\norm{u}.$$

\item[(4)] There is a constant $C>0$ such that for any $n\in \mathbb{N}$, $u : \ell^n_2 \rightarrow E$ and $v: E\rightarrow \ell^n_2$ we have $$\sup_{k\geq 1} ka_k(vu) \leq C\ell^*(v)\pi_{2,H}(u^*).$$

\end{itemize}
\end{thm}
\begin{proof}
(1) $\Rightarrow$ (2) : Consider $u : S \rightarrow H_n$ and $u^* : H_n \rightarrow E^*/S^\perp$. Since $E$ is weak type $(2,H)$,
it is $K$-convex as a Banach space, and so is $E^*$. Thus, by Lemma 11.7 of \cite{P1} there is $\widetilde{v} : H_n \rightarrow E^*$
such that $Q\widetilde{v} = u^*$ and $$\ell(\widetilde{v}) \leq 2K(E)\ell(u^*) \leq 2 \sqrt{n}\norm{u}K(E),$$
where $Q: E^* \rightarrow E^*/S^\perp$ is the natural quotient map.

Let $\widetilde{u} = \widetilde{v}^*|_E$. Then $\widetilde{u} : E\rightarrow H_n$ extends $u$ and
$$\ell(\widetilde{u}^*) = \ell(\widetilde{v}) \leq 2\sqrt{n}\norm{u}K(E).$$ Since $E^*$ is weak cotype $(2,H)$, we have
$$\sqrt{k}a^o_k(\widetilde{u}) \leq wC^H_2(E^*)\ell(\widetilde{u}^*)$$ for $k = 1, 2, \cdots, n$.
Then by the definition of cb-approximation number there is a projection $P_k$ on $H_n$ with $\text{rank} > n-k$ such that
$$\norm{P_k\widetilde{u}}_{cb} \leq wC^H_2(E^*)\ell(\widetilde{u}^*) \leq 2K(E)\sqrt{\frac{n}{k}}\norm{u},$$ which implies (2).

\vspace{0.3cm}

(2) $\Rightarrow$ (3) : We get (3) by setting $k=[(1-\delta)n]$, $C_{\delta} = C(1-\delta)^{-\frac{1}{2}}$, $P = P_k$ and $v = P_k \widetilde{u}$.

\vspace{0.3cm}

(3) $\Rightarrow$ (4) : Consider any $u : \ell^n_2 \rightarrow E$, $v : E \rightarrow \ell^n_2$, and fix $k \geq 1$.
Note that we can assume that $u$ is invertible by a typical perturbation argument. Now we set $F = \text{ran}(u) \subseteq E$.
Then by Theorem 5.8 of \cite{P1} there is a constant $C'>0$ such that $$\sqrt{k}c_k(v|_F) \leq C'\ell(i^*v^*) \leq C'\ell(v^*) \leq C'K(E)\ell^*(v),$$ where $i : F \hookrightarrow E$ and $c_k(\cdot)$ is the $k$-th Gelfand number defined by
$$c_k(T : X \rightarrow Y) = \inf \{\norm{T|_S} : S \subseteq X, \,\, \text{codim}(S) < k \}.$$
Thus we can choose $S \subseteq F$ such that $$\norm{v|_S} \leq 2C'K(E)k^{-\frac{1}{2}}\ell^*(v)\; \text{and}\; \text{dim}(F/S) < k.$$
By applying (3) to $v|_S$ we get a projection $$P : H_n \rightarrow H_n\; \text{with}\; \text{rk}P > \delta n$$ and an extension
$$w : E \rightarrow H_n\; \text{of}\; Pv|_S\; \text{with}\; \norm{w}_{cb} \leq C_{\delta}\norm{v|_S}.$$

Now we let $Q : H_n \rightarrow H_n$ be the orthogonal projection onto $u^{-1}(S)$.
Note that $$\text{rk}(I-Q) = n -\text{rk}Q \leq n - \text{dim}S < k$$ and $$\text{rk}(I-P) = n -\text{rk}P < (1-\delta)n.$$
Thus, for $\delta' = 1-\delta$ we have
\begin{align*}
a_{2k + [\delta' n]}(vu) & \leq a_{k + [\delta' n]}(Pvu) \leq a_{k}(PvuQ) = a_{k}(wuQ)\\
& \leq a_{k}(wu) = a_{k}(u^*w^*) = x^H_k(u^*w^*)\\
& \leq \norm{w^*}_{cb}x^H_k(u^*)\leq 2C_{\delta}C'K(E)k^{-1}\ell^*(v)\pi_{2,H}(u^*).
\end{align*}
The last inequality comes from Proposition \ref{Prop-Weyl}.

When $k = [(\frac{1-\delta'}{4})n]$, we have $$2k + [\delta' n] \leq [(\frac{1+\delta'}{2})n],$$ so that
\begin{align*}
\Big[\Big(\frac{1+\delta'}{2}\Big)n\Big]a_{[(\frac{1+\delta'}{2})n]}(vu)
& \leq \frac{4(1+\delta')}{1-\delta'}\Big[\Big(\frac{1-\delta'}{4}\Big)n\Big]a_{[(\frac{1+\delta'}{2})n]}(vu)\\
& \leq \frac{4(1+\delta')}{1-\delta'}ka_{2k+[\delta' n]}(vu) \leq C''\ell^*(v)\pi_{2,H}(u^*),
\end{align*}
where $C'' = \frac{8(1+\delta')}{1-\delta'}C_{\delta}C'K(E)$. Finally by Lemma \ref{upgrading-estimate-lemma2} we get our desired result.

\end{proof}

Before we establish equivalent formulations of weak type $(2,H)$ we need the following lemma.

\begin{lem}\label{lem-Lewis-variation}
Let $\text{dim}E = n$. Then there is an isomorphism $u: E \rightarrow H^*_n$ such that
$$\pi_{2,H}(u) = \sqrt{n}\; \text{and}\; \norm{u^{-1}}_{cb} = 1.$$
\end{lem}
\begin{proof}
By Lewis' theorem (Theorem 3.1. of \cite{P1}) there is an isomorphism $u: E \rightarrow H^*_n$ such that
$$\pi_{2,H}(u) = \sqrt{n}\; \text{and}\; \pi^*_{2,H}(u^{-1}) = \sqrt{n}.$$
Note that we have $$\pi^*_{2,H}(u^{-1}) = \inf\{\norm{B}_{HS}\norm{A}_{cb}\},$$
where the infimum runs over all possible factorization $u^{-1} : \ell^n_2 \stackrel{B}{\longrightarrow} H^*_n \stackrel{A}{\longrightarrow}E$.
Thus, we can apply the same argument as in the proof of Theorem 9.1 in \cite{P1.5}.
\end{proof}

Finally, we have the following equivalent formulation of weak cotype $(2,H)$.

\begin{thm}\label{thm-equi-weakOHtype2}
The followings are equivalent.

\begin{itemize}
\item[(1)] $E$ has weak type $(2,H)$.

\item[(2)] There are constants $0< \delta < 1$ and $C_{\delta}>0$ satisfying for any $S\subseteq E^{**}$ and any $u : S\rightarrow H_n$
there is a projection $P : H_n \rightarrow H_n$ of $\text{rank}>\delta n$ and $v : E^{**} \rightarrow H_n$ such that
$$v|_S = Pu\; \text{and}\; \norm{v}_{cb} \leq C_{\delta}\norm{u}.$$

\item[(3)] There is a constant $C>0$ such that for any $n\in \mathbb{N}$, $u : \ell^n_2 \rightarrow E^{**}$ and
$v : E^{**} \rightarrow \ell^n_2$ we have $$\sup_{k\geq 1}a_k(vu) \leq C \ell^*(v)\pi_{2,H}(u^*).$$

\end{itemize}
\end{thm}
\begin{proof}
(1) $\Rightarrow$ (2) : Since $E$ is $K$-convex, $E^{**}$ has weak type $(2,H)$ by duality (Proposition \ref{duality}).
Then we get our desired result by Theorem \ref{locweakOHtype2}.

\vspace{0.3cm}

(2) $\Rightarrow$ (3) : Theorem \ref{locweakOHtype2}.

\vspace{0.3cm}

(3) $\Rightarrow$ (1) : Note that $E$ is weak type 2 and consequently $K$-convex by Theorem 3.1 (c) \cite{P0}.
Now we consider $v : E \rightarrow H_n$ and choose $F \subseteq E^*$ with $\text{ran}v^* \subseteq F$ and $\text{dim}F = n.$
By Lemma \ref{lem-Lewis-variation} there is an isomorphism $$u : H^*_n \rightarrow F$$ such that
$$\norm{u}_{cb} = 1\;\text{and}\;\pi_{2,H}(u^{-1})\leq \sqrt{n}.$$
By Proposition 5.1. of \cite{P1.5} there is an extension $$\widetilde{u^{-1}} : E^{***} \rightarrow H_n$$ with
$$\pi_{2,H}(\widetilde{u^{-1}}) = \pi_{2,H}(u^{-1}).$$ If we let $w = \widetilde{u^{-1}}|_{E^*}$, then for any $k\geq 1$ we have
\begin{align*}
a^o_k(v) = a^o_k(v^*)& \leq \norm{u}_{cb}a^o_k(u^{-1}v^*) = a_k(wv^*)\\
& = a_k(v^{**}w^*) \leq Ck^{-1}\ell^*(v^{**})\pi_{2,H}(w^{**}).
\end{align*}
Note that $w^{**} = \widetilde{u^{-1}}$, so that we have
$$\sqrt{k}a^o_k(v) \leq C\sqrt{\frac{n}{k}}\ell^*(v^{**})\leq CK(E)\sqrt{\frac{n}{k}}\ell^*(v)$$ by applying (3) to $w^*$ and $v^{**}$.
If we set $k = [\frac{n}{2}]$, then we get our desired result by Lemma \ref{upgrading-estimate-lemma3}.

\end{proof}

\begin{rem}{\rm
We do not know whether we can replace $E^{**}$ into $E$ in Theorem \ref{thm-equi-weakOHtype2} at the time of this writing.
}
\end{rem}

\section{Weak-$H$ space}\label{sec-weak-H}

Now we consider an operator space analogue of weak Hilbert space and its various equivalent characterizations as in chapter 12. of
\cite{P1}. We assume that $H$ is perfectly Hilbertian, i.e. $H$ and $H^*$ are both subquadratic,
since we will consider weak type $(2,H)$ and weak cotype $(2,H^*)$ simultaneously.
Note that we use $d_{cb}(\cdot,\cdot)$ instead of $d_s(\cdot,\cdot)$ in (2) of Theorem \ref{thm-wOH-equi},
which looks more natural in the operator space setting.

\begin{defn}
$E$ is called a \textbf{weak-$H$ space} if it has both weak type $(2,H)$ and weak cotype $(2,H^*)$.
\end{defn}

\begin{rem}\label{wOH-reflexive}{\rm Clearly, a weak-$H$ space is weak Hilbert space as a Banach space and consequently reflexive
by Theorem 14.1 of \cite{P1}.}
\end{rem}

First, we present an analogue of Theorem 12.2. of \cite{P1} with a suitably modified proof.

\begin{thm}\label{thm-wOH-equi}
The followings are equivalent.
\begin{itemize}
\item[(1)] $E$ is a weak-$H$ space.

\item[(2)] For any $0<\delta <1$ there is a constant $C >0$ with the following property : for any finite dimensional $F\subseteq E$
we can find $F_1 \subseteq F$ and an onto projection $P :E \rightarrow F_1$ satisfying
$$d^H_{F_1,cb} := d_{cb}(F_1, H_{\text{dim}F_1}) \leq C, \; \text{dim}F_1 \geq \delta \text{dim}F\,\, \text{and} \,\, \norm{P}_{cb} \leq C.$$

\item[(3)] There are constants $0<\delta <1$ and $C >0$ with the same property as in (2).

\item[(4)] There are constants $0< \alpha < 1$ and $C>0$ such that for any $n\in \mathbb{N}$ and $u : H_n \rightarrow E$ we have
$$a^o_{[\alpha n]}(u) \leq C [\alpha n]^{-\frac{1}{2}}\pi_{2,H}(u^*).$$

\item[(5)] There is a constant $C>0$ such that for any $n\in \mathbb{N}$ and $u : H_n \rightarrow E$ satisfies
$$\sup_{k \geq 1} \sqrt{k}a^o_k(u) \leq C\pi_{2,H}(u^*).$$

\item[(6)] There is a constant $C>0$ such that for any finite dimensional $F\subseteq E^*$ with $\text{dim}F = n\in \mathbb{N}$
and for any $k = 1, \cdots, n$ we can find $F_1 \subseteq F$ with $\text{dim}F_1 > n-k$ and a projection $P : E^* \rightarrow F_1$
such that $$\gamma_{H^*}(P) \leq C\sqrt{\frac{n}{k}}.$$

\end{itemize}

Here, $\gamma_{H^*}(\cdot)$ is defined as follows. For $T\in CB(E,F)$ $$\gamma_{H^*}(T) := \inf \{ \norm{T_1}_{cb}\norm{T_2}_{cb}\},$$
where the infimum runs over all possible factorization $$T : E \stackrel{T_1}{\longrightarrow} H^*(I) \stackrel{T_2}{\longrightarrow} F$$
for some index set $I$.

\end{thm}
\begin{proof}
(1) $\Rightarrow$ (2) : First, we consider $C > 0$ and $C_{\delta}>0$ in (2) of Theorem \ref{locweakOHcotype2} and (3) of Theorem \ref{locweakOHtype2}.
Then for fixed finite dimensional subspace $F \subseteq E$ there is a further subspace $\widetilde{F} \subseteq F$
with $\text{dim}\widetilde{F} (=n)\geq \sqrt{\delta} \text{dim}F$ and an isomorphism $T : \widetilde{F} \rightarrow H_n$ with
$$\norm{T} \leq 1 \,\, \text{and} \,\, \norm{T^{-1}}_{cb}\leq C.$$
Now we apply (3) of Theorem \ref{locweakOHtype2} to $T$. Then there is a projection $Q : H_n \rightarrow H_n$ and a map $v : E \rightarrow H_n$
such that $$\text{rk}Q > \sqrt{\delta}\text{dim}\widetilde{F}, \, v|_{\widetilde{F}} = QT\;\text{and}\; \norm{v}_{cb} \leq C_{\delta}.$$
Set $F_1 := T^{-1}(\text{ran}Q)$ and $P := T^{-1}Qv$. Then clearly we have
$$d^H_{F_1,cb} \leq \norm{T^{-1}}_{cb}\norm{v}_{cb} \leq CC_{\delta}$$ and $\text{dim}F_1 = \text{rk}(Q) > \delta \text{dim}F$.
Moreover, for any $x\in OH_n$ we have $$PT^{-1}Qx = T^{-1}QvT^{-1}Qx = T^{-1}QQTT^{-1}Qx = T^{-1}Qx.$$
Thus, $P$ is a projection onto $F_1$ and $$\norm{P}_{cb} \leq \norm{T^{-1}}_{cb}\norm{v}_{cb} \leq CC_{\delta},$$
and consequently we get the desired result for $C' = CC_{\delta}.$

\vspace{0.3cm}

(2) $\Rightarrow$ (3) : Trivial.

\vspace{0.3cm}

(3) $\Rightarrow$ (4) : Let $\alpha = 1 -\frac{\delta}{2}$, and consider $u : H_n \rightarrow E$.
We may assume that $F=u(H_n)$ is of dimension $n$. By applying (3) to $F$ we get $F_1 \subseteq F$ and a projection $P : E \rightarrow F_1$ with $$\text{dim}F_1 = [\delta n] + 1, d^H_{F_1,cb} \leq C \; \text{and}\; \norm{P}_{cb} \leq C.$$
Then we can choose an isomorphism $T: F_1 \rightarrow H_m$ for $m = \text{dim}F_1$ with $$\norm{T}_{cb}\norm{T^{-1}}_{cb}\leq C.$$ Thus, we get
\begin{align*}
\Big( \sum_{k\geq 1}a^o_k(Pu)^2 \Big)^{\frac{1}{2}} & \leq \norm{T^{-1}}_{cb} \Big( \sum_{k\geq 1}a^o_k(TPu)^2 \Big)^{\frac{1}{2}}\\
& = \norm{T^{-1}}_{cb}\norm{TPu}_{HS} = \norm{T^{-1}}_{cb}\norm{u^*P^*T^*}_{HS}\\
& = \norm{T^{-1}}_{cb}\pi_{2,H}(u^*P^*T^*) \\ & \leq \norm{T^{-1}}_{cb}\norm{T}_{cb}\norm{P}_{cb}\pi_{2,H}(u^*)
\end{align*}
and consequently $$a^o_k(Pu) \leq C^2\pi_{2,H}(u^*)\sqrt{\frac{1}{k}}.$$ By choosing $k = \big[ \frac{\delta n}{2} \big]$ we get
$$a^o_{[\frac{\delta n}{2}]}(Pu) \leq C^2\pi_{2,H}(u^*)\Big[\frac{\delta n}{2}\Big]^{-\frac{1}{2}}.$$
Now we choose a map $v : H_n \rightarrow E$ such that $\text{rk}(v) = \big[ \frac{\delta n}{2} \big] - 1$ and
$$\norm{Pu-v}_{cb} \leq 2C^2\pi_{2,H}(u^*)\Big[\frac{\delta n}{2}\Big]^{-\frac{1}{2}}.$$
Since $\text{rk}(v+u-Pu) < n - \big[\frac{\delta n}{2}\big] -1$ and $u - (v+u-Pu) = Pu-v$, we get
$$a^o_{[\alpha n]}(u) \leq a^o_{n - [\frac{\delta n}{2}]}(u) \leq 2C^2\pi_{2,H}(u^*)\Big[\frac{\delta n}{2}\Big]^{-\frac{1}{2}}.$$
Note that $$\Big[\frac{\delta n}{2}\Big] \geq n - [\alpha n] \geq \Big(\frac{1}{\alpha} - 1\Big) [\alpha n].$$ Thus, we have
$$a^o_{[\alpha n]}(u) \leq 2 \Big(\frac{1}{\alpha} - 1\Big)^{-\frac{1}{2}}C^2\pi_{2,H}(u^*)[\alpha n]^{-\frac{1}{2}}.$$

\vspace{0.3cm}

(4) $\Rightarrow$ (5) : By Lemma \ref{upgrading-estimate-lemma}.

\vspace{0.3cm}

(5) $\Rightarrow$ (6) : Let $F\subseteq E^*$ be a subspace with $\text{dim}F =n$.
Then by Lemma \ref{lem-Lewis-variation} there is an isomorphism $T : H_n \rightarrow F$ such that
$$\norm{T}_{cb}=1\; \text{and}\;\pi_{2,H}(T^{-1})=\sqrt{n}.$$
By Proposition 5.1. of \cite{P1.5} there is a map $v : E^* \rightarrow H_n$ such that $$v|_F = T^{-1} \;\text{and}\;\pi_{2,H}(v) = \sqrt{n}.$$
Note that the condition (5) implies the condition (v) in Theorem 12.2 of \cite{P1}, which is equivalent to $E$ is a weak Hilbert space,
and consequently $E$ is reflexive by Theorem 14.1 of \cite{P1}. Thus $v = u^*$ for some $u : H_n \rightarrow E$ by the reflexivity of $E$.
Then we have by (5) that $$\sup_{k \geq 1} \sqrt{k}a^o_k(v) = \sup_{k \geq 1} \sqrt{k}a^o_k(u) \leq C\pi_{2,H}(u^*) = C\pi_{2,H}(v) = C\sqrt{n}.$$
From the definition of cb-approximation numbers there is a projection $Q : H_n \rightarrow H_n$ with $\text{rank} > n-k$ such that
$$\norm{Qv}_{cb} \leq C \sqrt{\frac{n}{k}}.$$ Indeed, there is a map $\widetilde{v} : E^* \rightarrow H_n$ such that
$$\text{rk}\widetilde{v} < k\; \text{and}\; \norm{v-\widetilde{v}}_{cb} \leq C \sqrt{\frac{n}{k}}.$$
Let $Q$ be the orthogonal projection from $H_n$ onto $\text{ran}\widetilde{v}^{\perp}$.
Then $Q$ is the desired projection since $$\norm{Qv}_{cb} = \norm{Qv -Q\widetilde{v}}_{cb} \leq C\sqrt{\frac{n}{k}}\;\,
\text{and}\;\, \text{rk}Q > n-k.$$

Now let $F_1 = \text{ran}TQ$. Then $\text{dim}F_1 > n-k$ and $P=TQv$ is a projection onto $F_1$. Moreover, we have
$$\gamma_{H^*}(P) \leq \norm{T}_{cb}\norm{Qv}_{cb}\leq C \sqrt{\frac{n}{k}}.$$

\vspace{0.3cm}

(6) $\Rightarrow$ (1) : First, we observe that (6) implies (2) for $(E^*, H^*)$. Indeed, if we let $P : E^* \rightarrow F_1$ be the projection in (6),
then there are maps $A: E^*\rightarrow H^*(I)$ and $B: H^*(I) \rightarrow F_1$ such that $$\norm{A}_{cb}\norm{B}_{cb} \leq C\sqrt{\frac{n}{k}}$$
for some index set $I$. Thus, if we set $k = [(1-\delta)n]$ when $(1-\delta)n \leq 1$ and $k = 1$ otherwise, we have
$$d^{H^*}_{F_1,cb} := d_{cb}(F_1, H^*_{\text{dim}F_1}) \leq C\sqrt{\frac{n}{k}}\leq C\sqrt{\frac{2}{1-\delta}},$$ which implies (2) for $E^*$.

Then by applying the implication (2)$\Rightarrow$(3)$\Rightarrow$(4)$\Rightarrow$(5)$\Rightarrow$(6) we get (6) for $(E^*, H^*)$,
and by applying the above result we get (2) for $(E^{**}, H)$. Now we have that $(E^{**}, H)$ and $(E^*, H^*)$ both satisfy condition (2),
which means that $E^*$ has weak type $(2,H^*)$ and weak cotype $(2,H)$ by Proposition \ref{duality}.

\end{proof}

We will consider an equivalent condition of weak-$H$ spaces using Grothendieck numbers corresponding (xi) of Theorem 12.6. in \cite{P1}.
Recall that $\Gamma_n(X)$, the $n$-th Grothendieck number of a Banach space $X$ is defined by
$$\Gamma_n(X) := \sup\Big\{\abs{\text{\rm det}(\left\langle x_i, y_j\right\rangle_{i,j})}^{\frac{1}{n}}
: (x_i)^n_{i=1} \subseteq B_X,\;(y_j)^n_{j=1} \subseteq B_{Y^*}\Big\},$$
where $B_X$ and $B_{X^*}$ imply corresponding unit balls of $X$ and $X^*$, respectively.
It is well known that a Banach space $X$ is a weak Hilbert space if an only if $$\sup_{n\geq 1} \Gamma_n(X) <\infty$$
((ix) of Theorem 12.6. of \cite{P1}).

There is a reformulation of $\Gamma_n(T)$ appropriate to consider an operator space analogue.
Let $$\Gamma'_n(X) := \sup\{\abs{\text{\rm det}(vu)}^{\frac{1}{n}}\},$$
where the supremum runs over any $u : \ell^n_2 \rightarrow X$ and $v: \ell^n_2 \rightarrow X$
with $$\pi_2(u^*)\leq \sqrt{n}\; \text{and}\; \pi_2(v) \leq \sqrt{n}.$$
By combining (3) of Proposition 1.1. and Theorem 2.1. in \cite{Ge} and Proposition 3.2. in \cite{DeJu} we have
$$\Gamma_n(X) \leq \Gamma'_n(X) \leq C \cdot \Gamma_n(X)$$ for some constant $C>0$ independent of $n$.

Now we define $\Gamma^H_n(E)$, the $n$-th $H$-Grothendieck number of an operator space $E$ by
$$\Gamma^H_n(E) := \sup\{\abs{\text{\rm det}(vu)}^{\frac{1}{n}}\},$$
where the supremum runs over any $u : \ell^n_2 \rightarrow E$ and $v: \ell^n_2 \rightarrow E$
with $$\pi_{2,H}(u^*)\leq \sqrt{n}\; \text{and}\; \pi_{2,H^*}(v) \leq \sqrt{n}.$$

We will also consider an equivalent condition of weak-$H$ spaces using eigenvalue estimations of certain linear maps
corresponding (ix) of Theorem 12.6. in \cite{P1}.
Let $\gamma^*_{H}$ be the trace dual of $\gamma_H$,
and note that we have for any finite rank linear map $T : E\rightarrow F$ between operator spaces that
$$\gamma^*_H(T) =\inf\{ \pi_{2,H^*}(A)\pi_{2,H}(B^*) \},$$
where the supremum runs over all possible factorizations $T : E \stackrel{A}{\longrightarrow} \ell_2(I) \stackrel{B}{\longrightarrow} F$
for some index set $I$ (Theorem 6.1. of \cite{P1.5}).

In the proof we need an analogue of Theorem 3.1 in \cite{P0}. Consequently, we add three more equivalent conditions of weak-$H$ spaces as follows.

\begin{thm}\label{thm-eigen-char}
The conditions (1), $\cdots$, (6) in Theorem \ref{thm-wOH-equi} are equivalent to the followings.
\begin{itemize}
\item[(7)] There is a constant $C>0$ such that for any $n\in \mathbb{N}$,
$u : \ell^n_2 \rightarrow E$ and $v: E \rightarrow \ell^n_2$ we have $$\sup_{k\geq 1} ka_{k}(vu) \leq C\pi_{2,H}(u^*)\pi_{2,H^*}(v).$$

\item[(8)] $$\sup_{n\geq 1} \Gamma^H_n(E) <\infty.$$

\item[(9)]
There is a constant $C>0$ such that $$\sup_{k\geq 1} k\abs{\lambda_k(T)} \leq C \cdot \gamma^*_H(T)$$
for all finite rank linear maps $T : E\rightarrow E$.
\end{itemize}
\end{thm}
\begin{proof}

(5) $\Rightarrow$ (7) : Consider $u : \ell^n_2 \rightarrow E$ and $v: E \rightarrow \ell^n_2$.
Then by Proposition \ref{Prop-Weyl-multiplicative} and \ref{Prop-Weyl} we have
\begin{align*}
\sup_{k\geq 1}ka_k(vu) & = \sup_{k\geq 1}kx^{H^*}_k(vu) \leq 2\sup_{n\geq 1}\sqrt{n}x^{H^*}_n(v)\sup_{k\geq 1}\sqrt{k}x^{H^*}_k(u)\\
& \leq 2 \pi_{2,H^*}(v)\sup_{k\geq 1}\sqrt{k}x^{H^*}_k(u) \leq 2 \pi_{2,H^*}(v)\sup_{k\geq 1}\sqrt{k}a^o_k(u : H_n \rightarrow E)\\
& \leq 2 \pi_{2,H^*}(v)\pi_{2,H}(u^*).
\end{align*}
\vspace{0.3cm}

(7) $\Rightarrow$ (5) : Consider $u : H_n \rightarrow E$ and choose $F \subseteq E$ with $\text{ran}u \subseteq F$ and $\text{dim}F = n$.
By Lemma \ref{lem-Lewis-variation} there is an isomorphism $v : F \rightarrow H_n$ such that
$$\pi_{2,H^*}(v) = \sqrt{n}\;\text{and}\;\norm{v^{-1}}_{cb} = 1.$$
Now we consider an extension $\widetilde{v} : E \rightarrow H_n$ of $v$ with $$\pi_{2,H^*}(\widetilde{v}) = \pi_{2,H}(v).$$ Then we have 
$$a^o_k(u) \leq \norm{v^{-1}}_{cb}a^o_k(\widetilde{v}u) = a_k(\widetilde{v}u) \leq Ck^{-1}\pi_{2,H}(u^*)\pi_{2,H^*}(\widetilde{v}) \leq C\frac{\sqrt{n}}{k}\pi_{2,H}(u^*).$$ Thus, for $k = [\frac{n}{2}]$ we get $$a^o_k(u)\leq 2\sqrt{2}Ck^{-\frac{1}{2}}\pi_{2,H}(u^*),$$
which leads us to the desired result by Lemma \ref{upgrading-estimate-lemma}.

\vspace{0.3cm}

(7) $\Rightarrow$ (8) : Consider any $u : \ell^n_2 \rightarrow E$ and $v: \ell^n_2 \rightarrow E$
with $$\pi_{2,H}(u^*)\leq \sqrt{n}\; \text{and}\; \pi_{2,H^*}(v) \leq \sqrt{n}.$$
Then by Proposition 1.3. of \cite{Ge} we have for some constant $C'>0$ that
\begin{align*}
\abs{\text{\rm det}(vu)}^{\frac{1}{n}} & = \Big(\prod^n_{k=1} a_k(vu)\Big)^{\frac{1}{n}}
\leq \Big(\frac{1}{n!}\Big)^{\frac{1}{n}} \sup_{k \geq 1}a_k(vu)\\
& \leq \frac{C'}{n}\sup_{k \geq 1}a_k(vu) \leq \frac{C'C}{n}\pi_{2,H}(u^*) \pi_{2,H^*}(v)\\
& \leq C'C.
\end{align*}

\vspace{0.3cm}

(8) $\Rightarrow$ (7) : Consider any $u : \ell^n_2 \rightarrow E$ and $v: \ell^n_2 \rightarrow E$. Then we have
\begin{align*}
\frac{n}{\pi_{2,H}(u^*) \pi_{2,H^*}(v)} a_n(vu) &\leq \frac{n}{\pi_{2,H}(u^*) \pi_{2,H^*}(v)}\Big(\prod^n_{k=1} a_k(vu)\Big)^{\frac{1}{n}}\\
& = \frac{n}{\pi_{2,H}(u^*) \pi_{2,H^*}(v)}\abs{\text{\rm det}(vu)}^{\frac{1}{n}}\\
&\leq \Gamma^H_n(E).
\end{align*}

\vspace{0.3cm}

(5) $\Rightarrow$ (9) : Let $T : E\rightarrow E$ be a finite rank linear map. For any given $\epsilon > 0$ we consider a factorization
$$T: E\stackrel{A}{\longrightarrow} \ell_2(I)\stackrel{B}{\longrightarrow} F$$ for some index set $I$ with
$$\pi_{2,H^*}(A)\pi_{2,H}(B^*) \leq (1+\epsilon)\gamma^*_H(T).$$ Then we have by \eqref{cb-approx-multiplicative}
$$a^o_{2k-1}(T)\leq a^o_k(B)a^o_k(A) = a^o_k(B^*)a^o_k(A) \leq C^2 k^{-1}\pi_{2,H^*}(A)\pi_{2,H}(B^*)$$
and consequently $$\sup_{k\geq 1}k\abs{\lambda_k(T)} \leq \sup_{k\geq 1}ka^o_k(T) \leq 4(1+\epsilon)\gamma^*_H(T).$$

\vspace{0.3cm}

(9) $\Rightarrow$ (8) : Consider any $u : \ell^n_2 \rightarrow E$ and $v: E \rightarrow \ell^n_2$ with
$$\pi_{2,H}(u^*)\leq \sqrt{n}\; \text{and}\; \pi_{2,H^*}(v) \leq \sqrt{n}.$$
Now we choose an unitary matrix $U$ such that $\abs{vu} = Uvu$. Then for some constant $C'>0$ we have
\begin{align*}
\abs{\text{\rm det}(vu)}^{\frac{1}{n}} & = \text{\rm det}(\abs{vu})^{\frac{1}{n}} = \Big[\prod^n_{k=1} \lambda_k(Uvu)\Big]^{\frac{1}{n}}\\
& = \Big[\prod^n_{k=1} \lambda_k(uUv)\Big]^{\frac{1}{n}} \leq C\gamma^*_H(uUv)\Big(\frac{1}{n!}\Big)^{\frac{1}{n}}\\
& \leq \frac{CC'}{n}\gamma^*_H(uUv) \leq \frac{CC'}{n}\pi_{2,H}(U^*u^*)\pi_{2,H^*}(v)\\ & \leq CC'.
\end{align*}

\end{proof}

\begin{rem}{\rm

\begin{itemize}
\item[(1)]
With the same proof as Theorem \ref{thm-eigen-char} we get the following equivalent formulations in the style of Theorem 3.1 in \cite{P0}.

\begin{itemize}
\item[$\bullet$]
$E$ has weak type $(2,H)$ if and only if there is a constant $C>0$ such that for any $n\in \mathbb{N}$,
$u : \ell^n_2 \rightarrow E^*$ and $v: E^* \rightarrow \ell^n_2$ we have $$\sup_{k\geq 1} ka_{k}(vu) \leq C\ell^*(u^*)\pi_{2,H^*}(v).$$

\item[$\bullet$]
$E$ has weak cotype $(2,H)$ if and only if there is a
constant $C>0$ such that for any $n\in \mathbb{N}$, $u : \ell^n_2 \rightarrow E$ and $v: E \rightarrow \ell^n_2$ we have
$$\sup_{k\geq 1} ka_{k}(vu) \leq C\ell(u)\pi_{2,H}(v).$$
\end{itemize}

\item[(2)]
Unlike in the Banach spaces cases the following condition is not equivalent to weak-$H$ spaces.
\begin{itemize}
\item[$(i)$]
There is a constant $C>0$ such that $$\sup_{k\geq 1} k\abs{\lambda_k(T)} \leq C \cdot \nu^o(T)$$
for all completely nuclear map $T : E\rightarrow E$, where $\nu^o(\cdot)$ means completely nuclear norm (\cite{ER}).
\end{itemize}
Indeed, if we consider another homogeneous Hilbertian operator space $\Hi$ which is not completely isomorphic to $H$,
then $\Hi$ satisfies $(i)$ but it is not a weak-$H$ space. Actually, $\Hi$ satisfies the following stronger condition. 

\begin{itemize}
\item[$(ii)$]
There is a constant $C>0$ such that $$\sum_{k\geq 1} \abs{\lambda_k(T)} \leq C \cdot \nu^o(T)$$
for all completely nuclear map $T : E\rightarrow E$.
\end{itemize}
This condition is obtained from the Banach space case and the fact that the set of all nuclear maps on $\Hi$ coincide with
the set of all completely nuclear maps on $\Hi$ with equivalent norms.

It is easy to see that the same proof works even though we replace the cotype $(2,H)$ condition into weak cotype $(2,H)$ in Proposition 3.8. of \cite{L}. Thus we can conclude that if $\Hi$ is a weak-$H$ space then it is completely isomorphic to $H$ itself, which is a contradiction.

\item[(3)]
If an infinite dimensional operator space $E$ satisfies the condition $(ii)$, then there is a constant $C'>0$ such that
$E$ is $C''$-homogeneous and $C''$-Hilbertian and $C'' \leq C'C^2$, i.e. $E$ is $C''$-isomorphic to a Hilbert space and every bounded linear map
$u : E\rightarrow E$ is completely bounded with $\norm{u}_{cb} \leq C''\norm{u}$.

Indeed, let $E$ be an infinite dimensional operator space satisfying the condition $(ii)$.
Then by Theorem 2.1. of \cite{P1.7} it is enough to show that
every finite dimensional subspace of $E$ is completely complemented with constants bounded by $C$.
Let $F$ be a finite dimensional subspace of $E$ and $j : F \hookrightarrow E$ be the inclusion.
Then by trace duality (Corollary 3.1.3.9 of \cite{Ju}) we have for any $S\in CB(F)$ that
\begin{align*}
\nu^o(jS) & = \sup\{ \abs{\text{tr}(vjS)} : \norm{v: E \rightarrow F}_{cb} \leq 1 \}\\
& = \sup\{ \abs{\text{tr}(uS)} : \norm{u}_{e,cb} \leq 1 \},
\end{align*}
where $$\norm{u}_{e,cb} := \inf\{ \norm{\tilde{u}}_{cb} : \tilde{u} \in CB(E,F), \tilde{u}|_{F} = u \}.$$
This implies $\norm{\cdot}_{e,cb}$ and $\nu^o(j\cdot)$ are norms on $B(F)$ which are in trace duality. Thus, we have
$$\norm{I_F}_{e,cb} = \sup\{ \abs{\text{tr}(S)} : S\in CB(F), \nu^o(jS) < 1 \}.$$
By Proposition 12.2.3 of \cite{ER} we have the following commutative diagram :
\begin{eqnarray*}
\begin{CD}
F @> jS >> E\\
@V A VV @AA B A\\
B(\ell_2) @>> M_{ab} > S_1
\end{CD}
\end{eqnarray*}
with $\norm{A}_{cb}\norm{B}_{cb} < \nu^o(jS)$ and $\norm{a}_{HS}, \norm{b}_{HS} < 1$,
where $M_{ab}$ is the two-sided multiplication operator by $a$ and $b$.
Since $B(\ell_2)$ is injective in the sense of operator space, we can extend $A$ to
$$\widetilde{A} : E \rightarrow B(\ell_2)\; \text{with}\; \norm{\widetilde{A}}_{cb} = \norm{A}_{cb}.$$
Now let $\widetilde{S} = BM_{ab}\widetilde{A} \in CB(E)$, then since $\widetilde{S}$ is an extension of $S$, eigenvalues of $S$
are also eigenvalues of $\widetilde{S}$. Thus, we have
\begin{align*}
\sum^{\infty}_{k=1} \abs{\lambda_k(S)} & \leq
\sum^{\infty}_{k=1} \abs{\lambda_k(\widetilde{S})} \leq
C\nu^o(\widetilde{S})\\ & \leq C
\norm{\widetilde{A}}_{cb}\norm{B}_{cb}\norm{a}_{HS}\norm{b}_{HS} < C,
\end{align*}
and consequently $$\norm{I_F}_{e,cb} < C.$$

\item[(4)]
By trace duality we have $$\gamma^*_H(T) \leq \nu^o(T)$$ for any finite rank map $T :E\rightarrow E$.
Thus, we can conclude that the following condition is equivalent to being completely isomorphic to $H$.
\begin{itemize}
\item[$(iii)$]
There is a constant $C>0$ such that $$\sum_{k\geq 1} \abs{\lambda_k(T)} \leq C \cdot \gamma^*_H(T)$$
for all completely nuclear map $T : E\rightarrow E$.
\end{itemize}
Indeed, $(iii)$ implies $E$ is $C''$-homogeneous and $C''$-Hilbertian for come constant $C'' > 0$ by (3)
and is a weak-$H$ space by Theorem \ref{thm-eigen-char}. By Proposition 3.8. of \cite{L} again (with a suitably modified proof) we have that
$E$ is completely isomorphic to $H$.

\end{itemize}

}
\end{rem}

\bibliographystyle{amsplain}
\providecommand{\bysame}{\leavevmode\hbox
to3em{\hrulefill}\thinspace}

\end{document}